\documentclass[letterpaper, 10 pt, conference]{ieeeconf}  
\IEEEoverridecommandlockouts  
\usepackage{graphicx} 
\usepackage[dvipsnames]{xcolor} 
\usepackage{amsmath} 
\usepackage{amssymb}  
\usepackage{amsthm,thmtools,dsfont} 
\usepackage{hyperref,cleveref} 
\let\labelindent\relax
\usepackage{enumitem} 
\usepackage{cite} 
\usepackage{microtype,enumerate,url} 
\usepackage{BOONDOX-ds} 
\usepackage{tikz} 
\usepackage[fulladjust]{marginnote} 
\usepackage{pgfplotstable} 
\usepackage{pgfplots} 
\usepackage{xstring} 
\usetikzlibrary{arrows,patterns} 
\usepackage{caption} 
\usepackage{mathabx} 
\usepackage{algorithm} 
\usepackage{algpseudocode} 
\usepgfplotslibrary{fillbetween} 
\usepackage{mathtools}
\definecolor{Gray}{gray}{0.9}
\definecolor{myRed}{RGB}{228,26,28}
\definecolor{myBlue}{RGB}{55,126,184}
\definecolor{myGreen}{RGB}{77,175,74}
\definecolor{myPurple}{RGB}{152,78,163}
\definecolor{myOrange}{RGB}{255,127,0}

\newcommand{\Tini}{{T_{\textup{ini}}}}
\newcommand{\Tf}{{T_{\textup{f}}}}
\newcommand{\wini}{{w_{\textup{ini}}}}
\newcommand{\wf}{w_{\textup{f}}}
\newcommand{\wref}{w_{\textup{ref}}}

\renewcommand{\theenumi}{(\roman{enumi})}

\makeatletter
\def\@IEEEsectpunct{.\ \,}
\def\paragraph{\@startsection{paragraph}{4}{\z@}{1.5ex plus 1.5ex minus 0.5ex}%
{0ex}{\normalfont\normalsize\itshape}}
\makeatother

\declaretheorem[style=definition]{theorem}

\declaretheorem[style=definition]{lemma}
\declaretheorem[style=definition,qed=$\vartriangle$]{remark}

\declaretheorem[style=definition]{problem} 

\declaretheorem[style=definition]{fact}

\newcommand {\nn}{\nonumber}
\newcommand{\beq}{\begin{equation}}
\newcommand{\eeq}{\end{equation}}
\newcommand {\bseq}{\begin{subequations}}
\newcommand {\eseq}{\end{subequations}}
\newcommand {\bma}{\left[}
\newcommand {\ema}{\right]}

\newcommand {\N}{\mathbb{N}} 	
\newcommand {\R}{\mathbb{R}} 	
\newcommand {\Rge}{\mathbb{R}_{+}} 	

\newcommand{\Ker}{\mathbf{ker} \,} 
\newcommand{\Image}{\mathbf{im} \,} 
\newcommand{\rank}{\mathbf{rank} \,} 
\newcommand{\transpose}{\mathsf{T}} 
\newcommand{\norm}[1]{\left\lVert#1\right\rVert}
\newcommand{\id}{\operatorname{id}}
\newcommand{\col}{\operatorname{col}}

\newcommand{\prox}{\mathbf{prox}} 
\newcommand{\zer}{\mathbf{zer}} 

\newcommand {\T}{\mathbb{T}} 	
\newcommand {\W}{\mathbb{W}} 	
\newcommand {\B}{\mathcal{B}} 	
\renewcommand{\L}{\mathcal{L}} 	
\newcommand{\pini }{{\pi_{\textup{ini}}}}
\newcommand{\pf}{{\pi_{\textup{f}}}}
\newcommand{\Pini }{{\Pi_{\textup{ini}}}}
\newcommand{\Pf}{{\Pi_{\textup{f}}}}

\usepackage[acronym]{glossaries}
\newacronym{Split-as-a-Pro}{Split-as-a-Pro}{operator \textbf{Split}ting \textbf{a}nd \textbf{s}calable \textbf{a}lternating \textbf{Pro}jections}
\newacronym{LQT}{LQT}{Linear Quadratic Tracking}
\newacronym{LTI}{LTI}{Linear Time-Invariant}
\newacronym{LQT-Split}{LQT-Split}{Linear Quadratic Tracking with Davis-Yin Splitting}
\newacronym{POCS}{POCS}{Projection Onto Convex Sets}
\newacronym{LQR}{LQR}{Linear Quadratic Regulator}
\newacronym{PID}{PID}{Proportional-Integral-Derivative}
\newacronym{DYS}{DYS}{DY Splitting}
\newacronym{QP}{QP}{Quadratic Programming}
\newacronym{DeePC}{DeePC}{Data-Enabled Predictive Control}
\newacronym{MPC}{MPC}{Model Predictive Control}
\newacronym{LQT-Split-as-a-Pro}{LQT-Split-as-a-Pro}{LQT with Split-as-a-Pro}
\newacronym{CCP}{CCP}{Convex, Closed, and Proper}
\newacronym{ITOS}{ITOS}{Inexact Three-Operator Splitting}
\newacronym{ADMM}{ADMM}{Alternating Direction Method of Multipliers}
\newacronym{SLS}{SLS}{System Level Synthesis}
\newacronym{LQT-FB-Split}{LQT-FB-Split}{LQT with FB Splitting}
\newacronym{LQT-DY-Split}{LQT-DY-Split}{LQT with DY Splitting}
\newacronym{FB}{FB}{Forward-Backward}
\newacronym{DY}{DY}{Davis-Yin}
\newacronym{LQT-FB-Split-as-a-Pro}{LQT-FB-Split-as-a-Pro}{LQT with FB Split-as-a-Pro}
\newacronym{LQT-DY-Split-as-a-Pro}{LQT-DY-Split-as-a-Pro}{LQT with DY Split-as-a-Pro}
\newacronym{SPC}{SPC}{Subspace Predictive Control}
\newacronym{DD-MPC}{DD-MPC}{Data-driven Model Predictive Control}

\title{\Large {\bfseries Split-as-a-Pro: behavioral control via operator splitting and alternating projections}}

\author{Yu Tang, Carlo Cenedese, Alessio Rimoldi, Florian D\"orfler, John Lygeros, Alberto Padoan\thanks{%
Y. Tang, C. Cenedese, A. Rimoldi, Florian D\"orfler, John Lygeros, are with the Department of Information Technology and Electrical Engineering at ETH Z\"urich (e-mails: {\tt\footnotesize  yutang@student.ethz.ch,\{ccenedese, arimoldi, dorfler, lygeros\}@control.ee.ethz.ch}). 
C. Cenedese is also with the Delft Center for Systems and Control at TU Delft. A. Padoan is with the Department of Electrical and Computer Engineering, University of British Columbia (e-mail: {\tt\footnotesize apadoan@ece.ubc.ca}). This work was supported as a part of NCCR Automation, a National Centre of Competence in Research, funded by the Swiss National Science Foundation (grant number 51NF40\_225155). 
  }
}

\date{\small\today} 

\begin{document}

\maketitle
\thispagestyle{empty}
\pagestyle{empty}

\begin{abstract}
The paper introduces \acrshort{Split-as-a-Pro}, a control framework that integrates behavioral systems theory, operator splitting methods, and alternating projection algorithms. The framework reduces dynamic optimization  problems — arising in both control and estimation —  to efficient projection computations.  \acrshort{Split-as-a-Pro} builds on a non-parametric formulation that exploits system structure to separate dynamic constraints imposed by individual subsystems from external ones — such as interconnection constraints and input/output constraints. This enables the use of arbitrary system representations, as long as the associated projection is efficiently computable, thereby enhancing scalability and compatibility with gray-box modeling.  We demonstrate the effectiveness of \acrshort{Split-as-a-Pro} by developing a  distributed  algorithm for solving finite-horizon linear quadratic control problems  and illustrate its use in predictive control.  Our numerical case studies show that algorithms obtained using \acrshort{Split-as-a-Pro} significantly outperform their centralized counterparts in runtime and scalability across various standard graph topologies, while seamlessly leveraging both model-based and data-driven system representations. 
\end{abstract}

\section{Introduction} 
Control across scales is one of the grand challenges in our field~\cite{annaswamy2023control}. Closing this gap is essential for controlling increasingly complex systems, like future energy grids, traffic networks and supply chains. While the literature on the subject is vast and opinions differ on the best approach,  there is broad consensus that effective solutions should exploit system structure, integrate prior knowledge when available, and combine model-based with data-driven control — particularly when some subsystems are well-modeled and others are not. 

The paper introduces \gls{Split-as-a-Pro}, a new control framework that merges non-parametric system modeling with advanced optimization tools,  offering a tractable and principled way to address complexity.  We base our approach on  behavioral systems theory~\cite{willems2007behavioral}, which models dynamical systems as sets of trajectories, independent of specific representations. 
The key idea is to reduce dynamic optimization problems to efficient \textit{projection} computations.  Leveraging operator splitting methods~\cite{bauschke2017convex,ryu2022large} and alternating projection algorithms~\cite{escalante2011alternating}, we decouple constraints imposed by individual subsystems  from external ones, such as interconnection and input/output constraints, allowing the projection to be computed efficiently. \linebreak   \gls{Split-as-a-Pro} offers several  advantages. First, non-parametric modeling based on behavioral system theory captures \emph{any} system representation as a special case and makes no distinction between inputs and outputs. This provides a rigorous foundation for networked control with gray-box models that leverage both model-based and data-driven representations.   Second, recognizing  projections  as the key bottleneck enables  a  systematic separation of constraints,  making it possible to exploit structure, parallelize computations, and design scalable control algorithms.  Third, \gls{Split-as-a-Pro} does not rely on a system being linear or finite-dimensional, nor on the cost function of being quadratic, thus allowing for broad generalizations to be explored. To illustrate our findings, we revisit the finite-horizon \gls{LQT} problem~\cite{anderson2007optimal}, a textbook example in control theory, yet one whose scalability still drives active research~\cite{wang2009fast,li2021frontiers}.

\subsubsection*{Related work}
Behavioral systems theory~\cite{willems2007behavioral} has gained renewed attention with the rise of data-driven control, largely due to Willems' fundamental lemma~\cite{willems2005note}. The lemma enables various data-driven control frameworks~\cite{depersis-formulas,henk-informativity,favoreel1999spc}
and, over finite time horizons, solutions to the \gls{LQT} problem~\cite{markovsky2008data}, which underpin recent direct predictive control algorithms~\cite{DeePC-2019,breschi-datadriven,berberich-datadriven}. 
The scalability of the \gls{LQT} problem has been widely studied for state-space systems~\cite{wang2009fast}, as it forms the backbone of \gls{MPC}~\cite{rawlings2009model}. Recent work has focused on enforcing locality constraints~\cite{li2021frontiers} through the \gls{SLS} parameterization~\cite{anderson2019system},  and on decoupling network behavior from individual subsystems using data-driven representations to enforce dissipativity properties~\cite{yan2021behavioural,yan2023distributed}. In contrast, we develop a scalable framework that supports both model-based and data-driven representations.

\subsubsection*{Contributions} The paper offers two main contributions. First, we introduce \gls{Split-as-a-Pro}, a new control framework that combines behavioral systems theory, operator splitting methods, and alternating projection algorithms.  Second, we present two algorithms that showcase the \gls{Split-as-a-Pro} framework allowing one to scale efficiently the solution of the \gls{LQT} problem to large systems, while seamlessly integrating model-based and data-driven representations, as well as handling constraints on inputs and outputs.

\subsubsection*{Paper organization} 
Section~\ref{sec:preliminaries} introduces notation, terminology, and preliminary results. Section~\ref{sec:problem-formulation} recalls the finite-horizon \gls{LQT} problem,  which serves as a motivating example throughout the paper. Section~\ref{sec:main-results} illustrates how the \gls{Split-as-a-Pro} framework solves the \gls{LQT} problem by distributing computations efficiently and incorporating additional constraints. Section~\ref{sec:simulations} presents numerical case studies that demonstrate our main results.   Section~\ref{sec:conclusion} concludes with a summary and an outlook on future research directions.

\section{Preliminaries} \label{sec:preliminaries}

\subsection{Notation and terminology}
The set of positive integer numbers is denoted by $\N$.   
The set of real and non-negative real numbers are denoted by $\R$ and $\Rge$, respectively. 
For ${T\in\N}$, the set of integers $\{1, 2, \dots , T\}$ is denoted by $\mathbf{T}$.  A map $f$ from $X$ to $Y$ is denoted by ${f:X \to Y}$; $(Y)^{X}$ denotes the set of all such maps. The \textit{restriction} of ${f:X \to Y}$ to a set ${X^{\prime}}$, with ${X^{\prime} \cap X \neq \emptyset}$, is denoted by $f|_{X^{\prime}}$ and is defined by $f|_{X^{\prime}}(x)$ for ${x \in X^{\prime}}$; if ${\mathcal{F} \subseteq (Y)^{X}}$, then $\mathcal{F}|_{X^{\prime}}$ is defined as ${\{ f|_{X^{\prime}} \, : \, f \in \mathcal{F}\}}$.
Notation for convex analysis is mostly borrowed from~\cite{ryu2022large}.
The gradient of a differentiable function $f$ is denoted by $\nabla f$. The subdifferential of a proper function $f$ is denoted by $\partial f$.  
The proximal operator of a \mbox{\gls{CCP}} function $f$ is denoted by $\prox_{f}$ and defined as in~\cite[Sec. 1.3.4]{ryu2022large}. The indicator function of the set $C$ is denoted by $\iota_{C}$ and the projection onto the set by $P_C$.
The identity operator is denoted by $\id$. The set of zeros of an operator $T$ is denoted by $\zer(T)$. 
 
\subsection{Behavioral systems theory}

\subsubsection{Sequences, shift operator, and Hankel matrices}
A \textit{sequence} is a function ${w: S \to \R^q}$, also denoted by $\{w_k\}_{k\in S}$, with ${S}$ a (non-empty) subset of consecutive elements of~$\N$; the sequence $w$ is \textit{finite} if $S$ is finite, \textit{infinite} otherwise. We use the terms \textit{sequence} and \textit{trajectory} interchangeably. By convention, a finite sequence ${w\in(\R^q)}^\mathbf{T}$ is often identified with the column vector ${w=\col(w(1),\ldots,w(T))\in\R^{qT}}.$ 

The \textit{concatenation} of sequences $w$ and $v$ is denoted by $w \wedge v$ and defined as in~\cite{markovsky2021behavioral}, with $w \wedge v = \col(w, v)$ if $w$ and $v$ are finite sequences identified with column vectors. 

The \textit{shift operator} ${\sigma : {(\R^q)}^\N \to {(\R^q)}^\N}$ is defined as $\sigma w (t)= w(t+1)$. By convention, the shift operator acts on all elements of a set ${\mathcal{W} \subseteq {(\R^q)}^\N}$, that is, $\sigma \mathcal{W}  =\{\sigma w\,:\, w\in\mathcal{W} \}$.

Given a permutation matrix ${\Pi\in\R^{q\times q}}$ and an integer \mbox{$0 < m < q$,} the map defined by the equation
\beq \label{eq:partition} 
(u,y) = \Pi w 
\eeq
induces a \textit{partition} of each ${w\in\R^q}$ into the variables ${u\in\R^m}$ and ${y\in\R^{q-m}}$.  We write $w \sim (u,y)$ if~\eqref{eq:partition} holds for some  permutation matrix ${\Pi\in\R^{q\times q}}$ and integer ${0 < m < q}$.  
Any partition of the form~\eqref{eq:partition} induces natural  projections  ${ \pi_u:  w \mapsto u}$ and ${ \pi_y:  w \mapsto y}$. 

The \textit{Hankel matrix} of depth ${L\in\mathbf{T}}$ associated with the sequence ${w\in {(\R^q)}^\mathbf{T}}$ is defined as 
\beq \label{eq:Hankel} 
\! H_{L}(w) \! = \!
\scalebox{1}{$
\bma  \nn
\begin{array}{ccccc}
w(1) & w(2)  & \cdots &  w(T-L+1)   \\
w(2) & w(3)  & \cdots &   w(T-L+2)   \\
\vdots  & \vdots  & \ddots & \vdots  \\
w(L) & w(L+1)  & \cdots  & w(T)
\end{array}
\ema
$} . \! \!
\eeq

\subsubsection{Systems}
In behavioral systems theory~\cite{willems2007behavioral}, a \textit{system} is a triple $\Sigma=(\T,\W,\B),$ where $\T$ is the \textit{time set}, $\W$ is the \textit{signal set}, and $\B \subseteq (\W)^{\T}$ is the \textit{behavior} of the system. Throughout this work, we exclusively focus on \textit{discrete-time} systems, with ${\T = \N}$ and ${\W = \R^q}$. We adapt definitions accordingly, emphasizing this assumption when necessary. We also identify each system $\Sigma$ with the corresponding behavior~$\B$.

\subsubsection{\gls{LTI} systems}
A system $\B$ is \textit{linear} if $\B$ is a linear subspace, \textit{time-invariant} if $\B$ is shift-invariant, \textit{i.e.}, ${\sigma(\B) \subseteq \B}$, and \textit{complete} if $\B$ is closed in the topology of pointwise convergence~\cite{willems1986timea}.
The set of all (discrete-time) complete \gls{LTI} systems is denoted by  $\L^q$. We often simply write $\B \in \L^{q}$.

\subsubsection{Kernel representations}
Every \gls{LTI} system ${\B \in \L^{q}}$ admits a \textit{kernel representation} of the form $\B = \Ker R(\sigma)$, where  $R(\sigma)$ is the operator defined by the polynomial matrix $R(z) = R_0 +R_1 z +\ldots+ R_{\ell}z^\ell,$ with $R_i\in\R^{p\times q}$ for $i\in\text{\boldmath$\ell$\unboldmath},$ and 
$\Ker R(\sigma) = \{w\in{(\R^q)}^\N \,:\, R(\sigma)w = 0\}.$ 
Without loss of generality, we assume  that  $\Ker R(\sigma)$  is a \emph{minimal}  kernel representation of $\B$,  \textit{i.e.},  $R(\sigma)$ has full row rank~\cite{willems1989models}.

\subsubsection{Integer invariants}
The structure of an LTI  system ${\B\in\L^q}$ is characterized by a set of \textit{integer invariants}~\cite{willems1986timea}, defined as
\begin{itemize}
\item the \textit{number of inputs} ${m(\B) = q-\text{row dim} R(\sigma)}$,  
\item the \textit{number of outputs} ${p(\B) = \text{row dim} R(\sigma)}$,  
\item the \textit{lag} ${\ell(\B) = \max_{i \in\mathbf{p}}\{\deg\text{row}_i R(\sigma) \}}$, and
\item the \textit{order} ${n(\B)=\sum_{i \in\mathbf{p}} \deg\text{row}_i R(\sigma) }$, 
\end{itemize}
where   
$\Ker R(\sigma)$  is a minimal  kernel representation of $\B$, while   
${\text{row dim}R}$ and ${\deg\text{row}_i R}$ are the number of rows and the degree of the $i$-th row of $R(z)$, respectively.  The integer invariants are intrinsic properties of a system, as they do not depend on its representation~\cite{willems1986timea}. Thus, we omit the dependence on the particular behavior if clear from context.

\subsubsection{State-space representations}
Every LTI system ${\B \in \L^{q}}$ can be described by the equations
\beq \label{eq:state-space}
\sigma x = Ax + Bu, \quad y=Cx+Du,
\eeq
and admits a (\textit{minimal}) \textit{input/state/output representation} 
\beq \label{eq:partition-ISO} 
\!  \! \B \!  = \!
\left\{  (u,y) \sim w \in (\R^{q})^\N \,:\, \exists \,x\in(\R^n)^{\N} \, \textup{s.t.}~\eqref{eq:state-space}~\text{holds}
\right\},
\eeq
where
$\scalebox{0.7}{$\bma\!
\begin{array}{cc}
A & B \\
C & D
\end{array}\!\ema $}
 \in \R^{(n+p)\times (n+m)}$ and  $m$, $n$, and $p$ are the number of inputs, the order, and the number of outputs, respectively. 

\subsubsection{Data-driven representations}
Data-driven representations describe system behaviors using raw data matrices, leveraging the fact that, over a finite time horizon, the behavior of an LTI system is a subspace with dimension defined by the integer invariants and time horizon. Next, we summarize a version of this principle known as the \textit{fundamental lemma}~\cite{willems2005note}.

\begin{lemma}~\cite{markovsky2021behavioral} \label{lemma:fundamental_generalized}
Let ${\B \in \L^{q}}$ and ${w \in \B|_{[1,T]}}$. Fix $L\in\mathbf{T}$, with ${L > \ell}$. Then  $\B|_{[1,L]} = \Image  H_L(w)$  if and only if
\beq  \label{eq:generalized_persistency_of_excitation} 
\rank  H_L(w)   =  mL+ n.
\eeq
\end{lemma}

\noindent
The rank condition~\eqref{eq:generalized_persistency_of_excitation} is referred to as the  \emph{generalized persistency of excitation} condition~\cite{markovsky2021behavioral}. Different variations of this principle can be formulated under a range of assumptions, see, \textit{e.g.}, the survey~\cite{markovsky2021behavioral} for an overview and~\cite{padoan2023data,berberich2023quantitative} for some recent extensions.

\section{Problem formulation}\label{sec:problem-formulation}

We consider the finite-horizon \gls{LQT} problem — a textbook example of control  design~\cite{anderson2007optimal}. Despite its simplicity, this example highlights the generality and tractability of our framework in addressing fundamental control  problems, the scalability of which is nontrivial and still an active area of research both in model-based and data-driven  contexts~\cite{wang2009fast,li2021frontiers}.%

Given an \gls{LTI} system ${\B\in\L^{q}}$,  the  goal   of \gls{LQT} is to find a trajectory ${\wf^{\star} \in\B}$ that is as close as possible to a given reference trajectory $\wref$.  
Over a finite-horizon, the problem can be formulated as that of minimizing the quadratic cost~\cite{markovsky2008data} 
\begin{equation}\label{eq:cost_quadratic}
\!\!
\norm{w-\wref}_{I \otimes \Phi}^2 \! = \! \sum_{t=1}^{\Tf} (w(t)-\wref(t))^{\transpose} \Phi (w(t)-\wref(t)) , \!
\end{equation}
where ${\Phi\in\R^{q\times q}}$ is a given symmetric positive definite matrix.

\begin{problem}[Finite-horizon \gls{LQT}]
Given an \gls{LTI} system $\B\in\L^{q}$, a reference trajectory $\wref\in\R^{q\Tf}$, an initial trajectory ${\wini \in \B|_{[1,\Tini]}}$, and a symmetric positive definite matrix ${\Phi\in\R^{q\times q}}$, find a trajectory ${\wf^{\star} \in \B|_{[1,\Tf]}}$ that minimizes the quadratic cost~\eqref{eq:cost_quadratic} and has $\wini$ as a prefix trajectory, \textit{i.e.},
\beq \label{eq:problem_linear_quadratic_tracking}
\begin{aligned}
\underset{\wf \in \R^{q\Tf}}{\min} \quad & \norm{\wf-\wref}_{I \otimes \Phi}^2 \\
\textup{s.t.} \quad & \wini \wedge \wf \in \B|_{[1,\Tini+\Tf]}.
\end{aligned}
\vspace{-.15cm}
\eeq
\end{problem}

\noindent 
For ${\Tini \ge \ell }$, the prefix trajectory $\wini$ implicitly fixes the initial condition for the optimal control problem~\eqref{eq:problem_linear_quadratic_tracking}, leading to a unique solution $\wf^{\star}$ that can be explicitly calculated~\cite{markovsky2008data}. This reduces to solving a constrained least squares problem, requiring $\Tf(m+p+n)(m+n)^2$  operations~\cite[p.368]{boyd2018introduction}. Complexity grows only linearly in the time horizon $\Tf$, but cubically in the order of the system  $n$   and the number of  inputs~$m$.   Moreover, this formulation does not readily integrate different system representations or exploit prior knowledge (e.g., interconnection topology). To address these limitations, we leverage the \gls{Split-as-a-Pro} framework, exploiting system structure via operator splitting methods~\cite{ryu2022large} and alternating projection algorithms~\cite{escalante2011alternating}.

\begin{remark}[Constraints] \label{remark:constraints}    
    While not part of the original \gls{LQT} problem, constraints are crucial for \gls{MPC} applications and recent data-driven variants~\cite{DeePC-2019,berberich-datadriven,breschi-datadriven}, which solve the finite-horizon \gls{LQT} problem in a receding-horizon fashion. Constraints can easily be incorporated into the \gls{LQT} problem~\eqref{eq:problem_linear_quadratic_tracking} by enforcing the additional constraint  $\wf \in \mathcal{C}$,  for a suitable choice of $\mathcal{C}$. As discussed in Section~\ref{sec:LQT-Split-AS-A-Pro}  and illustrated numerically in Section~\ref{sec:simulations},  constraints are readily handled by the \gls{Split-as-a-Pro} framework, provided that the projection onto $\mathcal{C}$ is ``simple,'' as is the case of box, non-negativity, or half-space constraints arising in predictive  control~\cite{parikh2014proximal}. 
\end{remark}

\section{Main results}\label{sec:main-results}

The \gls{Split-as-a-Pro} framework takes advantage of structural properties to efficiently solve dynamic optimization problems (e.g., control and filtering) using behavioral, non-parametric system representations. 
\gls{Split-as-a-Pro} involves three key steps summarized next and developed further individually in the remainder of the section.
\begin{itemize}
    \item[A.] \textit{Monotone inclusion problem}: Reformulate the original optimization problem as a monotone inclusion problem, possibly in a higher-dimensional space. 
    \item[B.] \textit{Operator splitting}: Select an operator splitting method to derive an iterative algorithm that solves the monotone inclusion problem, with convergence guarantees tied to the specific splitting method. 
    The main idea is to make the primary computational bottleneck the \textit{projection} onto the intersection of multiple closed, convex sets.
    \item[C.] \textit{Alternating projections}: 
    Decompose the projection task using an alternating projection algorithm into simpler projections onto individual sets, whose intersection forms the target set. This step reduces computational load by replacing a single complex projection with efficient, parallelizable projections,  isolating subsystem-specific constraints from interconnection or dynamic constraints.
\end{itemize}

\noindent  
Next, we apply this framework to the \gls{LQT} problem, showing how it enables the design of scalable algorithms that leverage non-parametric, finite-horizon system representations, compatible with both model-based and data-driven approaches.

\subsection{Monotone inclusion problem}
The first step in our analysis is to reformulate problem~\eqref{eq:problem_linear_quadratic_tracking} as a monotone inclusion problem.  To illustrate the flexibility of \gls{Split-as-a-Pro}, we present two possible reformulations.

\subsubsection{Monotone inclusion with two operators} 
The first is given by the (higher-dimensional) optimization problem
\beq \label{eq:2_objective}
\begin{aligned} 
\underset{w \in \R^{q(\Tini+\Tf)}}{\min} \  g(w) + h(w),
\end{aligned}
\eeq
where $g$ and $h$ are extended real-valued functions, defined as \vspace{-.5cm}
\bseq \label{eq:2_objective_splitting} 
\begin{align}
    g(w) &= \iota_{\B|_{[1,\Tini+\Tf]}}(w) + \iota_{\{\wini\}}(\pini(w)),\\
    h(w) &= \norm{\pf(w) - \wref}_{I \otimes \Phi}^2,
\end{align}
\eseq  
where the projections $\pini: w \mapsto \Pini w$ and $\pf: w \mapsto \Pf w$ are defined by the matrices $\Pini = [I_{q\Tini} \ 0_{q\Tini \times q\Tf}]$ and $\Pf = [0_{q\Tf \times q\Tini} \ I_{q\Tf}]$, respectively. 

This reformulation replaces the optimization variable $\wf$ with the higher-dimensional variable $w$, constrained to have a fixed prefix $\wini$. The constraints ${w \in \B|_{[1,\Tini+\Tf]}}$ and ${\pini(w) = \wini}$ are lifted into the objective via indicator functions. The corresponding monotone inclusion problem is
\beq \label{eq:2-operator-splitting}
0 \ \reflectbox{$\in$} \ \partial g + \nabla h .
\eeq

\begin{fact} \label{eq:1-1-correspondence-fg}
Consider the optimization problem~\eqref{eq:problem_linear_quadratic_tracking}
and the monotone inclusion problem~\eqref{eq:2-operator-splitting}, with $g$ and $h$ defined as in~\eqref{eq:2_objective_splitting}. Assume ${\Phi\in\R^{q\times q}}$ is symmetric and positive definite, ${\wini \in \B|_{[1,\Tini]}}$, and ${\Tini \ge \ell}$. Then $w^{\star} \in \zer(\partial g + \nabla h)$ if and only if $\pf(w^{\star})$  is a minimizer of~\eqref{eq:problem_linear_quadratic_tracking}.
\end{fact}

\noindent
Various reformulations can be tailored to use different operator splitting methods that take advantage of the particular properties of the terms appearing in the inclusion. 
 An alternative reformulation with three operators better suits scalability requirements and accommodates additional constraints, such as those arising in predictive control, as discussed next.

\subsubsection{Monotone inclusion with three operators}
Another possible equivalent reformulation of problem~\eqref{eq:problem_linear_quadratic_tracking} is given by the optimization problem
\beq \label{eq:3_objective}
\begin{aligned} 
\underset{w \in \R^{q(\Tini+\Tf)}}{\min} \ f(w) + g(w) + h(w) ,
\end{aligned}
\eeq
where $f$, $g$, and $h$ extended real-valued functions, defined as
\bseq \label{eq:DY_splitting_fgh}
\begin{align}
    f(w) &= \iota_{\B|_{[1,\Tini+\Tf]}}(w),\\
    g(w) &= \iota_{\{\wini\}}(\pini(w)),\\
    h(w) & = \norm{\pf(w) -\wref}_{I \otimes \Phi}^2.     
\end{align}
\eseq  
This reformulation mirrors the previous one,  but differs by explicitly separating the constraints ${w \in \B|_{[1,\Tini+\Tf]}}$ and ${\pini(w) = \wini}$, enabling distributed algorithms that scale efficiently and naturally accommodate additional constraints.    The corresponding monotone inclusion problem is
\beq \label{eq:3-operator-splitting}
0 \ \reflectbox{$\in$} \ \partial f + \partial g + \nabla h .
\eeq

\begin{fact} \label{eq:1-1-correspondence-fgh}
Consider the optimization problem~\eqref{eq:problem_linear_quadratic_tracking}
and the monotone inclusion problem~\eqref{eq:3-operator-splitting}, with $g$ and $h$ defined as in~\eqref{eq:DY_splitting_fgh}. Assume ${\Phi\in\R^{q\times q}}$ is symmetric and positive definite, ${\wini \in \B|_{[1,\Tini]}}$, and ${\Tini \ge \ell }$. Then $w^{\star} \in \zer(\partial f+ \partial g + \nabla h)$ if and only if $\pf(w^{\star})$  is a minimizer of~\eqref{eq:problem_linear_quadratic_tracking}. 
\end{fact}

\subsection{Operator splitting}

The second step in our framework 
is to solve the monotone inclusion problem through a fixed-point iteration of the form
\beq \label{eq:fixed-point-iteration}
w_{k+1} =  T(w_k) , \quad k \in \N,
\eeq 
with ${T}$ an operator defined by a given splitting method~\cite{ryu2022large}.

\subsubsection{\gls{FB} splitting}
One of the 
most well-known 
splitting methods is the Forward-Backward splitting~\cite[p.46]{ryu2022large}. 
 Since $g$ and $h$ are \gls{CCP} functions, with $h$ differentiable, this method  is applicable to problem~\eqref{eq:2_objective}. The solution is defined by the fixed-point iteration~\eqref{eq:fixed-point-iteration}, with $T$ defined as
\beq \label{eq:forward_backward}
\scalebox{1}{$
T= \prox_{\alpha \partial g}\left(\id - \alpha \nabla h \right).
$}
\eeq
Since $g$ is the indicator function of the affine set 
\beq \label{eq:Affine_Set} \!\!\!
\mathcal{A}|_{[1,\Tini+\Tf]} = \{  w \in \B|_{[1,\Tini+\Tf]} \,|\, \, w|_{[1,\Tini]} = \wini \} , \!
\eeq
this yields a standard \textit{projected gradient algorithm}~\cite[p.49]{ryu2022large}, described in pseudo-code in Algorithm~\ref{alg:LQT-FB-Split}.

\begin{algorithm}[h]
  \caption{}
  \label{alg:LQT-FB-Split}
  \textbf{Input:} Initial trajectory ${\wini \in \R^{q\Tini}}$,
  reference trajectory ${\wref \in \R^{q\Tf}}$,
  weight matrix ${\Phi \in \R^{q \times q}}$, 
  parameter ${\alpha \in \Rge}$, initial guess ${w_1 \in \R^{q(\Tf+\Tini)}}$. \\
  \textbf{Output:} Sequence of iterates $\{w_{k}\}_{k \in \N}$. 
  \begin{algorithmic}[1]
    \For{$k = 1, 2, \ldots$}
        \State \quad \label{alg:FB_LQT_Split_step1} $z_{k+1} = 2 \Pf^{\transpose} (I \otimes \Phi) (\Pf w_k -\wref)$ 
        \State \quad \label{alg:FB_LQT_Split_step2}
         $w_{k+1} = P_{\mathcal{A}|_{[1,\Tini+\Tf]}} \left(w_k -\alpha z_{k+1} \right)$
    \EndFor
  \end{algorithmic}
\end{algorithm}
\noindent
The next statement gives conditions for Algorithm~\ref{alg:LQT-FB-Split} to yield the desired solution. The proof is deferred to the appendix.

\begin{theorem}[Convergence of Algorithm~\ref{alg:LQT-FB-Split}] \label{thm:linear_quadratic_tracking_with_FB_splitting}
Consider the optimization problem~\eqref{eq:2_objective}, with $g$ and $h$ defined as in~\eqref{eq:2_objective_splitting}. Let ${\Phi \in \R^{q \times q}}$ be symmetric and positive definite, ${\wini \in \B|_{[1,\Tini]}}$, with ${\Tini \ge \ell}$, and ${\alpha \in (0,1/\rho(\Phi))}$. Then the sequence ${\{w_k\}_{k \in \N}}$ generated by Algorithm~\ref{alg:LQT-FB-Split} converges to the unique global optimum of the optimization problem~\eqref{eq:2_objective}.
\end{theorem}

\subsubsection{\gls{DY} Splitting}
Another particularly useful splitting for our purposes is the \gls{DY} splitting~\cite[p.48]{ryu2022large}.  Since   $f$, $g$, and $h$ are \gls{CCP} functions, with $h$ differentiable, this method applies to problem~\eqref{eq:3_objective}. The solution is defined by the fixed-point iteration \eqref{eq:fixed-point-iteration}, with $T$ defined as 
\beq \label{eq:davis_yin_operator}
\!\!
\scalebox{0.875}{$
T = \id - \prox_{\alpha \partial g} + \prox_{\alpha \partial f} \left(2 \prox_{\alpha \partial g} - \id - \alpha \nabla h \circ \prox_{\alpha \partial g}\right).
$} \!
\eeq
Since $f$ and $g$ are indicator functions of specific sets — the subspace $\B|_{[1,\Tini+\Tf]}$ and the singleton $\{\wini\}$ — the proximal operator of their subdifferentials reduces to a projection onto such sets~\cite[p.42]{ryu2022large}. The resulting iterative algorithm is described in pseudocode in Algorithm~\ref{alg:LQT-DY-Split}.

\begin{algorithm}[h]
  \caption{}
  \label{alg:LQT-DY-Split}
  \textbf{Input:} Initial trajectory ${\wini\in \R^{q\Tini}}$,
  reference trajectory ${\wref\in \R^{q\Tf}}$,
  weight matrix ${\Phi\in\R^{q\times q}}$, 
  parameter ${\alpha\in \Rge}$, initial guess ${w_1\in \R^{q(\Tf+\Tini)}}$. \\
  \textbf{Output:} Sequence of iterates $\{w_{k}\}_{k\in\N}$. 
  \begin{algorithmic}[1]
    \For{$k = 1, 2, \ldots, $}
        \State \quad \label{alg:LQT_Split_step1} $z_{k+1/2} = \wini \wedge \Pf w_k$
        \State \quad \label{alg:LQT_Split_step2} $z_{k+1} = 2 \Pf^{\transpose} (I \otimes \Phi) (\Pf z_{k+1/2} -\wref)$ 
        \State \quad \label{alg:LQT_Split_step3} $v_{k+1} ~= P_{\B|_{[1,\Tini+\Tf]}} \left(2z_{k+1/2} -w_k -\alpha z_{k+1}  \right)$
        \State \quad \label{alg:LQT_Split_step4} $w_{k+1} = w_{k} + v_{k+1} - z_{k+1/2}$
    \EndFor
  \end{algorithmic}
\end{algorithm}

\noindent
Next, we provide conditions for  Algorithm~\ref{alg:LQT-DY-Split}  to produce the desired solution. The proof can be found in the appendix.

\begin{theorem}[Convergence of Algorithm~\ref{alg:LQT-DY-Split}] \label{thm:linear_quadratic_tracking_with_davis_yin_splitting}
Consider the optimization problem~\eqref{eq:3_objective}, with $f$, $g$, and $h$ defined as in~\eqref{eq:DY_splitting_fgh}. Let ${\Phi\in\R^{q\times q}}$ be symmetric and positive definite, ${\wini \in \B|_{[1,\Tini]}}$, with ${\Tini \ge \ell}$, and ${\alpha \in (0,1/\rho(\Phi))}$. Then the sequence ${\{w_k\}_{k\in\N}}$ generated by Algorithm~\ref{alg:LQT-DY-Split} converges to the unique global optimum of the optimization problem~\eqref{eq:3_objective}.
\end{theorem}

\noindent
By design, Algorithms~\ref{alg:LQT-FB-Split} and~\ref{alg:LQT-DY-Split} yield the desired solution to the associated problems. All steps are highly parallelizable and conducive to scalable implementation, particularly when the cost is separable (e.g., $\Phi$ block-diagonal)~\cite{parikh2014proximal}. The only exception is the projection onto the sets 
$\mathcal{A}|_{[1,\Tini+\Tf]}$ and $\mathcal{B}|_{[1,\Tini+\Tf]}$, which poses a fundamental scalability challenge. For an \gls{LTI} system ${\B \in \L^{q}}$, the projection can be efficiently computed either in closed form or by solving a linear system in saddle point form~\cite{benzi2005numerical}. If ${B \in \R^{q(\Tini+\Tf) \times r}}$ is a matrix with orthonormal columns spanning the subspace $\mathcal{B}|_{[1,\Tini+\Tf]}$, the projection is given by  
\[
P_{\mathcal{B}|_{[1,\Tini+\Tf]}}(w) = B(B^{\transpose}B)^{-1} B^{\transpose}w.
\]
Thus, computing the projection $P_{\mathcal{B}|_{[1,\Tini+\Tf]}}$ requires inverting a $r \times r$ matrix and $2r^3/3$ operations~\cite[p.210]{boyd2018introduction}.  Since 
$r = \dim \mathcal{B}|_{[1,\Tini+\Tf]}  = m(\Tini+\Tf) + n$~\cite{willems2005note}, the complexity grows cubically with the number of inputs $m$ and the system order $n$, posing a challenge in predictive control applications, where projections must be computed repeatedly in a receding-horizon fashion. In such settings, the projection matrix can be precomputed offline, reducing each projection to a matrix-vector multiplication costing $q^2(\Tini+\Tf)^2$ operations~\cite[p.4]{golub2013matrix}. Yet, the offline computation can be prohibitively expensive—or even impossible—when dealing with large-scale systems with subsystems that are accessible only through local input-output data. Similar complexity considerations apply to the projection onto the affine set $\mathcal{A}|_{[1,\Tini+\Tf]}$. To address these limitations, we use alternating projections to exploit system structure, distribute projection computations, and support the simultaneous use of model-based and data-driven system representations, as discussed next.

\subsection{Alternating projections} \label{sec:LQT-Split-AS-A-Pro}

The third step in our analysis leverages alternating projection algorithms to exploit system structure.

\subsubsection{Alternating projection algorithms}
Alternating projection algorithms aim to find a point in the intersection of two sets~\cite{escalante2011alternating}.   A classical alternating projection algorithm, introduced by von Neumann~\cite{escalante2011alternating}, is defined by the iteration
\beq \label{eq:von-Neumann-iteration}
w_{k+1}=P_{\mathcal{C}_2}\left(P_{\mathcal{C}_1}\left(w_{k}\right)\right),  \quad k \in \N,
\eeq
where $\mathcal{C}_1$ and $\mathcal{C}_2$ are closed subspaces of a Hilbert space with non-empty intersection.  The main advantage is that, while projecting onto the intersection may be costly, projections onto individual sets are often efficient. This idea generalizes to more than two sets. A well-known  extension for multiple sets is Dykstra's algorithm~\cite[p.556]{bauschke2017convex}, defined by the iteration
\beq \label{eq:dykstra-iteration}
w_{k+1}=P_{\mathcal{C}_s}\left(P_{\mathcal{C}_{s-1}}\left(\cdots P_{\mathcal{C}_1}\left(w_k\right)\cdots\right)\right), \quad k \in \N.
\eeq
where $\mathcal{C}_1, \mathcal{C}_2, \dots, \mathcal{C}_s$, 
are closed convex subsets of a Hilbert space with non-empty intersection.

Alternating projection algorithms can be used to replace the monolithic projection steps in Algorithms~\ref{alg:LQT-FB-Split} and~\ref{alg:LQT-DY-Split}. The key idea is to decompose the sets $\mathcal{A}|_{[1, \Tini + \Tf]}$ and $\mathcal{B}|_{[1, \Tini + \Tf]}$ into sets onto which projections can be computed efficiently.
For example, this applies when $\mathcal{B}|_{[1, \Tini + \Tf]}$ arises from an interconnected \gls{LTI} system, which one can model as the intersection of the set of all isolated subsystem behaviors $\mathcal{C}_1$ with the set $\mathcal{C}_2$ of all interconnection constraints~\cite{yan2021behavioural,yan2023distributed}. The decomposition simplifies the projection step and lends itself to parallel implementation, which is advantageous in large-scale settings. In Algorithm~\ref{alg:LQT-FB-Split}, this replaces the projection onto the affine set~\eqref{eq:Affine_Set} with projections onto three sets: two sets $\mathcal{C}_1$ and $\mathcal{C}_2$ whose intersection define $\mathcal{B}|_{[1, \Tini + \Tf]}$ together with the affine set defined by the constraint  ${\pini(w) = \wini}$. In Algorithm~\ref{alg:LQT-DY-Split}, the projection onto $\mathcal{B}|_{[1, \Tini + \Tf]}$ is similarly replaced by alternating projections onto $\mathcal{C}_1$ and $\mathcal{C}_2$.  

The alternating projection framework also extends naturally to constraints not originally part of the \gls{LQT} problem. For instance, additional constraints can be incorporated by introducing convex sets $\mathcal{C}_3, \dots, \mathcal{C}_s$ and applying Dykstra's iteration~\eqref{eq:dykstra-iteration} over $\mathcal{C}_1, \mathcal{C}_2, \dots, \mathcal{C}_s$. As noted in Remark~\ref{remark:constraints}, such constraints are standard in \gls{MPC} and recent data-driven variants, and often include box or non-negativity constraints (e.g., for saturation or physical limits). These constraints fit naturally within the \gls{Split-as-a-Pro} framework, as the required projections are straightforward to compute,  e.g., for standard sets like box, non-negativity, or half-space constraints~\cite[Section 6]{parikh2014proximal}.  Section~\ref{sec:simulations} provides numerical illustrations of this extension in the context of data-driven predictive control. 

\subsubsection{\gls{LQT} with \gls{Split-as-a-Pro}}
Returning to the \gls{LQT} problem, we now illustrate how to replace single-step projections with alternating projection algorithms.  For illustration, we consider Algorithm~\ref{alg:LQT-DY-Split} and substitute the projection in Step~\ref{alg:LQT_Split_step3} with von Neumann's iteration~\eqref{eq:von-Neumann-iteration}. This yields a distributed version of Algorithm~\ref{alg:LQT-DY-Split},  whose pseudocode is given in Algorithm~\ref{alg:LQT-Split-as-A-Pro}.

\begin{algorithm}[th!]
  \caption{} 
  \label{alg:LQT-Split-as-A-Pro}
\textbf{Input:}  Initial trajectory ${\wini\in \R^{q\Tini}}$,
  reference trajectory ${\wref\in \R^{q\Tf}}$,
  weight matrix ${\Phi\in\R^{q\times q}}$, 
  parameters ${\alpha\in \Rge}$, number of inner loop iterations ${J \in\N}$, initial guess ${w_1\in \R^{q(\Tf+\Tini)}}$. \\
  \textbf{Output:} Sequence of iterates $\{w_{k}\}_{k\in\N}$. 
  \begin{algorithmic}[1]
    \For{$k = 1, 2, \ldots, $}
        \State \quad $z_{k+1/2} =   \wini \wedge \Pf w_k$ \label{alg3:step2}
        \State \quad $z_{k+1} = \Pf^{\transpose} (I \otimes \Phi) (\Pf z_{k+1/2} -\wref)$ 
        \State \quad $v_{k+1,0}=2z_{k+1} -w_k -2\alpha z_{k+1/2}$
        \For{$j =  0, 1, 2, \ldots, J-1 $} \label{alg:LQT_Split_as_a_Pro_Inner_Loop_step0} 
            \State \quad \label{alg:LQT_Split_as_a_Pro_Inner_Loop_step1} $v_{k+1,j+1/2} = P_{\mathcal{C}_1} \left(v_{k+1,j}  \right)$
            \State \quad \label{alg:LQT_Split_as_a_Pro_Inner_Loop_step2} $v_{k+1,j+1} = P_{\mathcal{C}_2} \left(v_{k+1,j+1/2}  \right)$
            
        \EndFor \label{alg:LQT_Split_as_a_Pro_Inner_Loop_step3} 
        \State \quad $w_{k+1} = w_{k} + v_{k+1,J} - z_{k+1/2}$
    \EndFor
  \end{algorithmic}
\end{algorithm}

\noindent
The inner loop of Algorithm~\ref{alg:LQT-Split-as-A-Pro} performs $J$ iterations, each involving two projections. If projections onto the individual sets are efficient, with worst-case per-iteration cost $O(r)$, the total cost scales as $O(Jr)$. This can offer significant advantages when the interconnection structure is sparse and projections decompose well. A detailed analysis of how to exploit such structure and implement the associated projections efficiently is deferred to future work. As previously discussed, the alternating projection step can accommodate additional constraints by introducing further sets $\mathcal{C}_3, \mathcal{C}_4, \dots, \mathcal{C}_n$ and applying Dykstra’s iteration~\eqref{eq:dykstra-iteration}. Section~\ref{sec:simulations} presents empirical evidence that Algorithm~\ref{alg:LQT-Split-as-A-Pro} outperforms centralized methods in sparse interconnection settings and, despite its simplicity, yields control performance on par with state-of-the-art data-driven predictive control algorithms.

\begin{remark}[Early termination]
Residual errors due to
early termination of the inner loop are inevitable in Algorithm~\ref{alg:LQT-Split-as-A-Pro}. Thus, it is crucial to determine whether Algorithm~\ref{alg:LQT-DY-Split} converges after a minimum number of alternating projection iterations while accounting for these errors. With a judicious choice of inner loop iterations $J$, convergence can be established using Fej\'er monotonicity~\cite{bauschke2017convex} or inexact operator splitting methods~\cite{ryu2022large}.  The detailed proof of convergence is deferred to future work due to space constraints. 
\end{remark}

\section{Numerical case studies} \label{sec:simulations}

We now illustrate the scalability of \gls{Split-as-a-Pro} by comparing the (centralized) Algorithm~\ref{alg:LQT-DY-Split} and the (distributed) Algorithm~\ref{alg:LQT-Split-as-A-Pro} on three common graph topologies — chain, ring, and lattice — by increasing the number of subsystems $\nu$ and analyzing runtime performance on a fixed \gls{LQT} problem. 

\subsubsection{Scalable control for interconnected system}
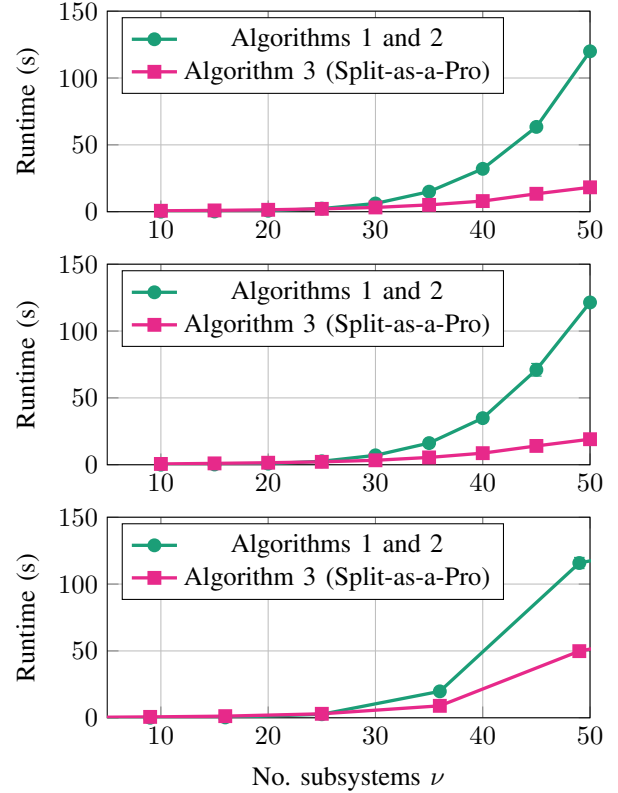
\begin{figure}[h!]
\centering

\begin{tikzpicture}
    \begin{axis}[
        ylabel={Runtime (s)},
        width=8cm, %
        height=4.25cm, 
        xmin=5, xmax=50, 
        ymin=0, ymax=150, 
        y label style={yshift=-5pt}, 
        grid=major,
        legend pos=north west,
    ]
    
    \addplot[
        color={rgb,255:red,27;green,158;blue,119},
        mark=*,
        line width=1.25pt, 
        error bars/.cd,
        y dir=both, y explicit,
    ] 
    table [
        x=units, 
        y=centralized_lqr_mean, 
        y error=centralized_lqr_var, 
        col sep=comma
    ] {10-50ecc,chain.csv};
    \addlegendentry{Algorithms~\ref{alg:LQT-FB-Split} and~\ref{alg:LQT-DY-Split}}

    \addplot[
        color={rgb,255:red,231;green,41;blue,138},
        mark=square*,
        line width=1.25pt, 
        error bars/.cd,
        y dir=both, y explicit,
    ] 
    table [
        x=units, 
        y=distributed_worst_mean, 
        y error=distributed_worst_var, 
        col sep=comma
    ] {10-50ecc,chain.csv};
    \addlegendentry{Algorithm~\ref{alg:LQT-Split-as-A-Pro} (\gls{Split-as-a-Pro})}
    \end{axis}
\end{tikzpicture}

\begin{tikzpicture}
    \begin{axis}[
        ylabel={Runtime (s)},
        width=8cm, %
        height=4.25cm, 
        xmin=5, xmax=50, 
        ymin=0, ymax=150, 
        y label style={yshift=-5pt}, 
        grid=major,
        legend pos=north west,
    ]

    
    \addplot[
        color={rgb,255:red,27;green,158;blue,119},
        mark=*,
        line width=1.25pt, 
        error bars/.cd,
        y dir=both, y explicit,
    ] 
    table [
        x=units, 
        y=centralized_lqr_mean, 
        y error=centralized_lqr_var, 
        col sep=comma
    ] {10-50ecc,ring.csv};
    \addlegendentry{Algorithms~\ref{alg:LQT-FB-Split} and~\ref{alg:LQT-DY-Split}}

    \addplot[
        color={rgb,255:red,231;green,41;blue,138},
        mark=square*,
        line width=1.25pt, 
        error bars/.cd,
        y dir=both, y explicit,
    ] 
    table [
        x=units, 
        y=distributed_worst_mean, 
        y error=distributed_worst_var, 
        col sep=comma
    ] {10-50ecc,ring.csv};
    \addlegendentry{Algorithm~\ref{alg:LQT-Split-as-A-Pro} (\gls{Split-as-a-Pro})}
    \end{axis}
\end{tikzpicture}

\begin{tikzpicture}
    \begin{axis}[
        xlabel={No. subsystems $\nu$},
        ylabel={Runtime (s)},
        width=8cm, %
        height=4.25cm, 
        xmin=5, xmax=50, 
        ymin=0, ymax=150, 
        y label style={yshift=-5pt}, 
        grid=major,
        legend pos=north west,
    ]

    
    \addplot[
        color={rgb,255:red,27;green,158;blue,119},
        mark=*,
        line width=1.25pt, 
        error bars/.cd,
        y dir=both, y explicit,
    ] 
    table [
        x=units, 
        y=centralized_lqr_mean, 
        y error=centralized_lqr_var, 
        col sep=comma
    ] {10-50ecc,grid.csv};
    \addlegendentry{Algorithms~\ref{alg:LQT-FB-Split} and~\ref{alg:LQT-DY-Split}}

    \addplot[
        color={rgb,255:red,231;green,41;blue,138},
        mark=square*,
        line width=1.25pt, 
        error bars/.cd,
        y dir=both, y explicit,
    ] 
    table [
        x=units, 
        y=distributed_worst_mean, 
        y error=distributed_worst_var, 
        col sep=comma
    ] {10-50ecc,grid.csv};
    \addlegendentry{Algorithm~\ref{alg:LQT-Split-as-A-Pro} (\gls{Split-as-a-Pro})}
    \end{axis}
\end{tikzpicture}
\centering
\caption{Runtime of  Algorithms~\ref{alg:LQT-FB-Split},~\ref{alg:LQT-DY-Split} and~\ref{alg:LQT-Split-as-A-Pro}  as a function of the number of interconnected subsystems $\nu$ for chain (top), ring  (center), and lattice (bottom) topologies.}
\label{fig:lqr-comparison}
\end{figure}

We consider an interconnected system inspired by~\cite{anderson2019system}, which consists of ${\nu\in\N}$ subsystems described by the state-space representations
\beq \label{eq:chain_dynamics}
\!\!\!
\begin{aligned}
    x_i(t+1) &= A_{i} x_i(t) + B_i u_i(t) + \sum_{j \not = i} B_i u_{i,j}(t), \\[-0.25cm]
y_i(t) &= C_i x_i(t) ,
\end{aligned} \!
\eeq
where ${x_i(t)\in\R^2}$,  ${u_i(t)\in\R}$, ${u_{i,j}(t)\in\R}$, ${y_i(t)\in\R}$,
and 
\begin{align*}
A_{i} &=
\scalebox{0.85}{$
\bma
\begin{array}{cc}
    1 & \Delta t \\
    -\tfrac{K_i}{m_i} \Delta t & 1 -\tfrac{d_i}{m_i} \Delta t 
\end{array}
\ema $}, \, 
B_i =
\scalebox{0.85}{$
\bma
\begin{array}{cc}
    0 \\
    1
\end{array}
\ema $}, \, 
C_i =
\scalebox{0.85}{$
\bma
\begin{array}{cc}
    \tfrac{k_{i}}{m_i} \Delta t & 0
\end{array}
\ema $}, 
\end{align*}
with ${m_i>0}$, ${d_i>0}$, ${k_i>0}$, ${K_{i} >0}$, and ${\Delta t>0}$. 
To illustrate compatibility with gray-box modeling, we model half of the subsystems using state-space representations, while the other half are represented by Hankel matrices derived from suitably collected data~\cite{markovsky2021behavioral}. The parameters of each subsystem are sampled from a uniform distribution. We select  \( \Tini = 2\nu+1 \) and $\Tf=5$. We collect input-output data \( w_{\textup{d}}  \sim(u_{\textup{d}}, y_{\textup{d}}) \) by applying pseudo-random inputs from a uniform distribution in $[-1, 1]$,  with \( T = \nu(\Tini+\Tf)+200 \) samples. 
The parameters in  Algorithms~\ref{alg:LQT-FB-Split},~\ref{alg:LQT-DY-Split} and~\ref{alg:LQT-Split-as-A-Pro}  are selected as $\alpha=0.1$ and $\Phi = I$, respectively. For Algorithm \ref{alg:LQT-Split-as-A-Pro}, we set $J=5$, which empirically generates satisfactory controller performance. Figure \ref{fig:lqr-comparison} compares the mean runtime of  Algorithms~\ref{alg:LQT-FB-Split},~\ref{alg:LQT-DY-Split} and~\ref{alg:LQT-Split-as-A-Pro}  across the three graph topologies averaged over five different runs. For all three graph topologies,  Algorithms~\ref{alg:LQT-FB-Split} and~\ref{alg:LQT-DY-Split}  are outperformed by  Algorithm~\ref{alg:LQT-Split-as-A-Pro}   even for moderately large-scale systems.

\subsubsection{Beyond linear quadratic tracking} \label{sec:control-perform}
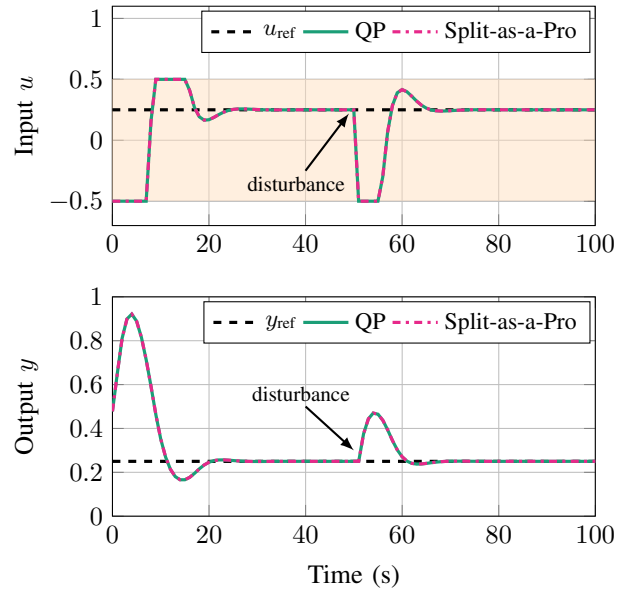
\begin{figure}[b!]
\centering
\begin{tikzpicture}
    \begin{axis}[
        width=0.45\textwidth,
        grid=both,
        minor ytick={-0.5,0.5},
        ylabel={Input $u$},
        legend style={
            font=\small,         
            at={(0.99,0.98)},     
            legend columns=3     
        },
        height=4.5cm, 
        xmin=0, xmax=100, 
        ymin=-0.7, ymax=1.1, 
        y label style={yshift=-2.5pt}, 
    ]
    \addplot[fill=myOrange!25, opacity = 0.5, draw=none, forget plot] coordinates {(0, -0.5) (100, -0.5) (100, 0.5) (0, 0.5)} -- cycle;
    \addplot[very thick, dashed, color=black] table[x=time, y=u_ref, col sep=comma] {w_ecc.csv};     
    \addlegendentry{${u_{\textup{ref}}}$}
    
    \addplot[very thick, solid, color={rgb,255:red,27;green,158;blue,119}] table[x=time, y=u_method_1, col sep=comma] {w_ecc.csv};
    \addlegendentry{QP}
    \addplot[very thick, dash dot, color={rgb,255:red,231;green,41;blue,138}] table[x=time, y=u_method_2, col sep=comma] {w_ecc.csv};
    \addlegendentry{\gls{Split-as-a-Pro}}
    \draw[-latex, thick, black] (axis cs:40,-0.25) -- (axis cs:49,0.22) node[below, text=black, xshift=-20pt, yshift=-20pt] {\footnotesize disturbance};
    \end{axis}
\end{tikzpicture}\\[.25cm]
\begin{tikzpicture}
    \begin{axis}[
        width=0.45\textwidth,
        grid=both,
        xlabel={Time (s)},
        ylabel={Output $y$},
                legend style={
            font=\small,         
            at={(0.99,0.98)},     
            legend columns=3     
        },
        height=4.5cm, 
        xmin=0, xmax=100, 
        ymin=-0., ymax=1., 
        y label style={yshift=-2.5pt}, 
    ]
    \addplot[very thick, dashed, color=black] table[x=time, y=y_ref, col sep=comma] {w_ecc.csv};
    \addlegendentry{$y_{\textup{ref}}$}
    
    \addplot[very thick, solid, color={rgb,255:red,27;green,158;blue,119}] table[x=time, y=y_method_1, col sep=comma] {w_ecc.csv};
    \addlegendentry{{QP}}
    \addplot[very thick, dash dot, color={rgb,255:red,231;green,41;blue,138}] table[x=time, y=y_method_2, col sep=comma] {w_ecc.csv};
    \addlegendentry{\gls{Split-as-a-Pro}}
    \draw[-latex, thick, black] (axis cs:40,0.50) -- (axis cs:50,0.30) node[above, text=black, xshift=-20pt, yshift=15pt] {\footnotesize disturbance};
    \end{axis}
\end{tikzpicture}
\caption{Control performance of \gls{DD-MPC} using a \gls{QP} solver (solid) and \gls{Split-as-a-Pro} (dashdotted), together with reference values (dashed) and input constraints set (shaded).}
\label{fig:control-perform-box}
\end{figure}

The \gls{Split-as-a-Pro} framework extends well beyond \gls{LQT} problems  — for instance, to input-constrained problems.  For illustration, we use a variant of Algorithm~\ref{alg:LQT-Split-as-A-Pro}, which includes an additional projection step onto the input constraint set, to solve the \gls{DD-MPC} problem from~\cite{berberich-datadriven} in a receding horizon fashion.  We consider the \gls{LTI} system described above with $\nu=2$ subsystems arranged in a chain graph topology. 
The reference is selected as ${\wref \sim (u_{\textup{ref}}, y_{\textup{ref}})},$
 with ${u_{\textup{ref}} = 0.25}$ and ${y_{\textup{ref}} = 0.25}$. The cost is defined by the block diagonal matrix $\Phi=\textup{diag}(0.1I,10I)
$. 
We add input box constraints \( u_i(t)\in[-0.5,\, 0.5] \).  Thus,   we leverage  von Neumanns's alternating projection algorithm~\eqref{eq:von-Neumann-iteration},  with \( \mathcal{C}_1 = \mathcal{B}|_{[1,\Tini+\Tf]} \) and \( {\mathcal{C}_2 = [-0.5, 0.5]^{m(\Tini+\Tf)} }\). We also add a disturbance to the system at the time ${t=50}$. Fig.~\ref{fig:control-perform-box} shows the input/output trajectory of the first subsystem obtained using the \gls{DD-MPC} algorithm given in~\cite{berberich-datadriven} with a \gls{QP} solver (solid) and \gls{Split-as-a-Pro} (dashdotted), together with reference values (dashed) and constraints (shaded). The results show that \gls{Split-as-a-Pro} attains control performance on par with those obtained using a QP solver, while effectively enforcing input constraints and adapting to exogenous disturbances.

\clearpage
\section{Conclusion} \label{sec:conclusion}
The paper has introduced \gls{Split-as-a-Pro}, a control framework that integrates behavioral systems theory with operator splitting methods and alternating projection algorithms. The framework reduces large-scale control problems to the efficient computation of a projection, offering advantages in terms of scalability and compatibility with gray box modeling.  We have showcased the benefits of \acrshort{Split-as-a-Pro} by developing a distributed algorithm for solving finite-horizon linear quadratic control problems, showing that it significantly outperforms its centralized counterparts for large-scale systems across various graph topologies and match state-of-the-art data-driven predictive control algorithms.  Future research will extend the \acrshort{Split-as-a-Pro} framework  to a wider class of dynamic optimization problems.

\bibliographystyle{IEEEtran}
\bibliography{refs}

\appendix
\subsection{Proofs of Theorems~\ref{thm:linear_quadratic_tracking_with_FB_splitting} and ~\ref{thm:linear_quadratic_tracking_with_davis_yin_splitting}}
We recall sufficient conditions for convergence of the fixed point iteration \eqref{eq:fixed-point-iteration} with $T$ defined as in~\eqref{eq:davis_yin_operator}~\cite[Ch. 2]{ryu2022large}. 
\begin{lemma}[Convergence of \gls{DY} splitting algorithm] \label{lemma:davis_yin}
Consider the optimization problem~\eqref{eq:3_objective}.  Assume 
\begin{itemize}
\item[(i)] $f$ and $g$ are \gls{CCP} functions,
\item[(ii)] $h$ is a convex function with $\gamma$-Lipschitz gradient,
\item[(iii)] $\zer(\partial f + \partial g +\nabla h) \not =  \emptyset$,
\item[(iv)] $\alpha \in (0,2/\gamma)$.
\end{itemize}
Then the sequence $\{w_{k}\}_{k\in\N}$ defined by~\eqref{eq:fixed-point-iteration}, with $T$ defined as in~\eqref{eq:davis_yin_operator}, converges to some ${w^{\star} \in \zer(\partial f + \partial g +\nabla h)}$. 
\end{lemma}
\noindent
If $g=0$ in \eqref{eq:davis_yin_operator} then we retrieve \eqref{eq:forward_backward} where $\partial g\coloneqq \partial f$. This lemma also implies convergence of the sequence attained from  ~\eqref{eq:fixed-point-iteration}, with $T$ as in \eqref{eq:forward_backward}, to ${w^{\star} \in \zer(\partial g +\nabla h)}$ under the same assumptions on $g$ and $h$. We now prove Theorem~\ref{thm:linear_quadratic_tracking_with_davis_yin_splitting}; the proof of Theorem~\ref{thm:linear_quadratic_tracking_with_FB_splitting} is analogous and, hence, omitted.\\
\textit{Proof of Theorem~\ref{thm:linear_quadratic_tracking_with_davis_yin_splitting}: }
First, we show that the assumptions of Lemma~\ref{lemma:davis_yin} hold, with $f$, $g$, and $h$ defined as in~\eqref{eq:DY_splitting_fgh}.

(i) 
If a set $C \subseteq \mathbb{R}^n$ is convex, closed, and non-empty, then its indicator function $\iota_C$ is \gls{CCP}~\cite[p.8]{ryu2022large}. Thus, 
since $f=\iota_{\B|_{[1,\Tini+\Tf]}}$ and $\B|_{[1,\Tini+\Tf]}$ is a linear subspace of $\R^{q(\Tini+\Tf)}$~\cite{markovsky2021behavioral}, $f$ is \gls{CCP}~\cite[p.8]{ryu2022large}. 
Similarly,  $g$ is convex, closed, since its epigraph is closed, and proper, as $g$ is finite if $w = \col(\wini, \wf)$ for some $\wf$ and equals ${+\infty}$, otherwise.

(ii) Satisfied by design since $\Phi$ is positive definite.

(iii) Since ${\wini \in \B|_{[1,\Tini]}}$ and ${\Tini \ge  \ell}$, any trajectory of the form ${\wini \wedge \wf  \in \B|_{[1,\Tini+\Tf]}}$ is feasible, ensuring
 $\zer(\partial f + \partial g +\nabla h) \neq \emptyset$.

(iv) This condition is satisfied by design.

\noindent
By Lemma~\ref{lemma:davis_yin}, $w_k\to w^{\star}\in \zer(\partial f + \partial g +\nabla h)$ as $k \to \infty$. 
Uniqueness of ${w^{\star}}$ follows from $\Phi$ being positive definite and the one-to-one correspondence established by Fact~\ref{eq:1-1-correspondence-fgh} between solutions of \eqref{eq:problem_linear_quadratic_tracking} and~\eqref{eq:3_objective}.\hfill\QED

\end{document}